\documentclass[12pt,reqno,english,a4]{amsart}
\usepackage{amsmath,amsthm,amsfonts,amssymb,amscd}
\usepackage[latin1]{inputenc}
\usepackage{psfrag}
\usepackage{epsfig}
\usepackage{a4wide}

\numberwithin{equation}{section}

\newcommand{\leb}{\operatorname{Leb}}
\newcommand{\dist}{\operatorname{dist}}

\newcommand{\Chi}{{\bf 1}}

\newcommand{\al} {\alpha}       
\newcommand{\be} {\beta}        
       
\newcommand{\de} {\delta}       

\newcommand{\e} {\epsilon}
\newcommand{\vare}{\varepsilon}

\newcommand{\f}{\varphi}

\def \RR {{\mathbb R}}
\def \ZZ {{\mathbb Z}}
\def \NN {{\mathbb N}}

\newcommand{\dem}{\begin{proof}}
\newcommand{\cqd}{\end{proof}}

\newcommand{\qand}{\quad\text{and}\quad}

\newcommand{\cc}{{\mathcal C}}
\newcommand{\cf}{{\mathcal U}}
\newcommand{\cv}{{\mathcal V}}
\newcommand{\cn}{{\mathcal N}}

\newcommand{\cp}{{\mathcal P}}

\newtheorem{maintheorem}{Theorem}

\newcommand{\cmt}{\begin{maintheorem}}
\newcommand{\fmt}{\end{maintheorem}}

\newtheorem{maincorollary}[maintheorem]{Corollary}

\newcommand{\cmc}{\begin{maincorollary}}
\newcommand{\fmc}{\end{maincorollary}}

\newtheorem{T}{Theorem}[section]
\newcommand{\cte}{\begin{T}}
\newcommand{\fte}{\end{T}}

\newtheorem{Corollary}[T]{Corollary}
\newcommand{\cco}{\begin{Corollary}}
\newcommand{\fco}{\end{Corollary}}

\newtheorem{Proposition}[T]{Proposition}
\newcommand{\cpr}{\begin{Proposition}}
\newcommand{\fpr}{\end{Proposition}}

\newtheorem{Lemma}[T]{Lemma}
\newcommand{\cle}{\begin{Lemma}}
\newcommand{\fle}{\end{Lemma}}

\newcommand{\csle}{\begin{Lemma}}
\newcommand{\fsle}{\end{Lemma}}

\theoremstyle{remark}

\newtheorem{Remark}[T]{Remark}
\newcommand{\cre}{\begin{Remark}}
\newcommand{\fre}{\end{Remark}}

\newtheorem*{Definition}{Definition}
\newcommand{\cde}{\begin{Definition}}
\newcommand{\fde}{\end{Definition}}

\begin{document}

\author{José F. Alves}
\address{Departamento de Matem\'atica Pura, Faculdade de Ci\^encias da Universidade do
Porto, Rua do Campo Alegre 687, 4169-007 Porto, Portugal} \email{jfalves@fc.up.pt}

\author{Krerley Oliveira}
\address{Departamento de Matemática, Universidade Federal de Alagoas,
Campus A. C. Simões s/n, 57000-000 Maceió, Brazil} \email{krerley@mat.ufal.br}

\author{Ali Tahzibi}
\address{Departamento de Matem\'atica, ICMC-USP São Carlos, Caixa
Postal 668, 13560-970 São Carlos, Brazil.} \email{tahzibi@icmc.sc.usp.br}

\subjclass{37C40, 37C75, 37D25}

\keywords{SRB measures, entropy, induced maps, non-uniform expansion}

\date{\today}

\thanks{Work carried out at USP-São Carlos, IMPA and University of Porto.
JFA was partially supported by FCT through CMUP, A.T was supported by FCT on leave from USP-São Carlos and K.O by Fapeal/Brazil and Pronex-CNPq/Brazil.}

\setcounter{tocdepth}{1}

\begin{abstract}
We consider classes of dynamical systems admitting Markov induced
maps. Under general assumptions, which in particular guarantee the
existence of SRB measures, we prove that the entropy of the SRB
measure varies continuously with the dynamics. We apply our result
to a vast class of non-uniformly expanding maps of a compact
manifold and prove the continuity of the entropy of the SRB
measure. In particular, we show that the SRB entropy of Viana maps
varies continuously with the map.

\end{abstract}

\title[On the continuity of the SRB entropy]{On the continuity of the SRB entropy\\ for endomorphisms}

\maketitle

\tableofcontents


\section{Introduction}

In this work we address ourselves to the study of the continuity of the metric entropy for endomorphisms.
Entropy of dynamical systems can be regarded quite generally as a measure of unpredictability. Topological
entropy measures the complexity of a dynamical system in terms of the exponential growth rate of the number of
orbits which can be distinguished over long time intervals, within a fixed small precision. Kolmogrov-Sinai's
metric entropy is an invariant which, roughly speaking, measures the complexity of the dynamical system in
probabilistic terms with respect to a fixed invariant measure.

Considering that the observable properties are, in a physical sense, the properties which hold on a positive
volume measure set, one tries to verify the existence of invariant measures with ``good" densities with respect
to the volume measure. Let us explain this in more precise terms.  We consider discrete-time systems, namely,
iterates of smooth transformations $f : M \rightarrow M$ on a Riemannian manifold. We consider a probability
measure  defined by a volume form on $M$ that we call {\em Lebesgue measure}. A Borel probability measure $\mu$
on $M$ is said to be a {\em Sinai-Ruelle-Bowen (SRB) measure} or a {\em physical measure}, if there exists a
positive Lebesgue measure subset of points $x \in M$ for which
 \begin{equation}\label{SRB}
 \lim_{ n
\rightarrow \infty } {\frac{ 1 }{n }} \sum_{ j = 0 }^{ n - 1 } \varphi ( f^j ( x ) ) = \int \varphi d \mu ,\quad
\text{for every }\varphi \in C^0 ( M ).
 \end{equation}
 The set of points $x\in M$ for which \eqref{SRB} holds is called the {\em basin} of the SRB measure~$\mu$.
Finding SRB measures  for a given dynamical system  may be  a difficult task  in general. By Birkhoff's ergodic
theorem, one possible way to prove the existence of these physically relevant measures is to construct
absolutely continuous invariant ergodic probabilities.  This kind of measures is  constructed in \cite{ABV} for
a vast class of diffeomorphisms and endomorphisms satisfying some weak hyperbolicity conditions.

Recently, there is an increasing emphasis on the study of the stability of the statistical properties of
dynamical systems. One natural formulation for this kind of stability corresponds to the continuous variation of
the SRB measures. Another interesting question in this direction is to ask whether the entropy of the SRB
measure varies continuously as a function of the dynamical system. The question of the continuity of the entropy
(topological or metric) is an old issue, going back to the work of Newhouse \cite{N}, for example.

It is known that uniformly expanding $C^2$ maps of a compact manifold admit a unique SRB measure which is
absolutely continuous with respect to Lebesgue measure and its density  varies continuously in the $L^1$ norm.
By means of this continuity and the entropy formula for these systems one  easily obtains the continuity of the
SRB entropy. For Axiom~A diffeomorphisms the continuity of SRB measures and  even more regularity is established
in \cite{R2} and  \cite{Man90}. The regularity of the SRB entropy for Axiom A flows is proved in \cite{C}.
Analiticity of metric entropy for Anosov diffeomorphisms is proved in \cite{Po}.

In this paper we present an abstract model and give sufficient conditions which imply the continuous variation
of the SRB entropy in quite general families of maps, including maps with critical sets.  Under the same
hypotheses, the continuous variation of the SRB measures  is proved in \cite{Al2}. It is important to remark
that in the presence of critical points it is not clear whether the continuous variation of absolutely
continuous invariant measures implies the continuous variation of their entropy or not. Let us observe that if
we do not have  absolute continuity,  the continuous variation of the SRB measures does not imply the continuity
of their entropy. For instance, in the quadratic family $f_a(x)= 4ax(1-x)$ one can find parameters $a$ for which
$f_a$ has an absolutely continuous SRB measure, and there is a sequence $a_n$ converging to $a$ with $f_{a_n}$
having a unique SRB measure concentrated on an attracting periodic orbit (sink). Furthermore, the Dirac measures
supported on those sinks converge to the SRB measure of $f_a$. This shows that the convergence of SRB measures
does not necessarily imply the convergence of the SRB entropy.

In the sequel we show that a large class  of {\em non-uniformly expanding} endomorphisms (admitting critical
sets) satisfy the conditions of our main result. We just suppose some natural \emph{slow recurrence} to the
critical set to construct the absolutely continuous invariant measures as in \cite{ABV}. We apply our results to
an open set of non-uniformly expanding endomorphisms constructed by Viana \cite{V}, and prove the continuity of
the entropy of the unique absolutely continuous invariant measure for such endomorphisms.

As far as we know our result is the first one  giving continuity of the SRB entropy for families of
endomorphisms admitting critical points. Our approach is different from the usual ways to prove the continuity
of the entropy. We construct induced maps for endomorphisms and relate the entropy of the SRB measure of the
initial system and the entropy of a corresponding measure of the induced system. Then we prove some continuity
results for the induced map and come back to the original map.

\subsubsection*{Acknowledgements} We are thankful to Marcelo Viana for several valuable discussions on these topics.

\section{Statement of results}

Let $M$ be a $d$-dimensional compact Riemannian manifold and denote the Lebesgue measure on $M$ by   $m$. We are
interested in studying the continuity of the metric entropy of smooth maps $f\colon M\to M$  with respect to
some physically relevant measure on $M$.

A very important tool that we will be using are \emph{induced maps}. Roughly speaking, an induced map for a
system $f$ is a transformation $F$ from some region of the ambient space into itself, defined for each point as
an iterate of $f$, where the number of iterations depends on the point. If we carry out this process carefully,
some asymptotic properties of $f$ (asymptotic expansion, for instance)  can be verified as properties of $F$ at
the first iteration (real expansion) for almost all points. A hard problem is to decode back the information
obtained for $F$ into information about the original dynamical system.

\subsection{Induced maps}

Let $ F:\Delta\rightarrow \Delta$ be an {\em induced map}\index{induced map} for $f$ defined in some topological
disk $\Delta\subset M$, meanning that there exists a countable partition $\cp$ of a full Lebesgue measure subset
of $\Delta$, and there exists a {\em return time} function $\tau\colon\cp\rightarrow\ZZ^+$ such that
 $$
 F\vert_\omega=f^{\tau(\omega)}\vert_\omega, \quad\mbox{ for each }\quad \omega\in
\cp.
 $$
We assume that the following conditions on the induced map $F$ hold:
 \begin{enumerate}
 \item[(i$_1$)]   Markov: {\em \(F\vert_\omega: \omega \to \Delta\)
 is a \( C^{2} \) diffeomorphism, for each \( \omega\in\mathcal P
\).} \item[(i$_2$)]  Uniform expansion: {\em there exists \(0<\kappa<1 \) such that
    for any \( \omega\in\cp \) and $x\in \omega$
    \[ \|DF(x)^{-1}\|<\kappa. \]}
\item[(i$_3$)]  Bounded distortion:
   {\em  there exists \( K >0 \) such that for any
    $\omega\in\cp$ and
    \( x,y\in \omega \) 
    \[
   \left|\frac{\det DF(x)}{\det DF(y)}-1\right| \leq
   K \dist(F(x),F(y)).
    \]}
\end{enumerate}

It is well known  that a map $F$ in these conditions has a unique absolutely continuous ergodic invariant
probability measure. Moreover, such a probability measure is equivalent to the Lebesgue measure on $\Delta$, and
its density is bounded from above and from below by constants. Proofs of these assertions will be given in
Proposition~\ref{l.density}. In this setting, we also prove in Proposition \ref{pr.entropyF} that if
 $F\colon\Delta\to\Delta$ is a piecewise expanding Markov induced map and $\mu_F$
is its absolutely continuous invariant probability measure,  then the entropy of $F$ with respect to the
probability measure $\mu_F$ satisfies:
 \begin{equation}\label{Theorem A}h_{\mu_F}(F)=\int_\Delta \log|\det DF(x)|\,d\mu_F.
 \end{equation}


A natural question is how to obtain an absolutely continuous $f$-invariant probability measure from the
existence of such measure for $F$. The integrability of the return time function $\tau\colon\Delta\to\ZZ^+$ with
respect to the Lebesgue measure on $\Delta$ is enough for the existence of this measure. Indeed, if $\mu_F$ is
the absolutely continuous $F$-invariant probability measure, then
\begin{equation}\label{eq.relmu}
 \mu_f^*
=\sum_{j=0}^{\infty}f_{\ast}^j\left(\mu_F\mid \{\tau_f>j\}\right)
 \end{equation}
is an absolutely continuous $f$-invariant finite  measure. We denote by $\mu_f$ the probability measure which is
obtained from $\mu_f^*$ by dividing it by its mass. Throughtly this paper we are assuming the integrability of the return time.

A formula similar to the  one displayed in \eqref{Theorem A} holds for $C^2$ endomorphisms $f$ of a compact manifold~$M$
with respect to an absolutely continuous invariant probability measure $\mu_f$. In fact, by \cite[Remark 1.2]{L}
the Jacobian function $\log |\det Df(x)| $ is always integrable with respect to  $\mu_f$. Then, by \cite[Theorem
1.1]{QZ}, if
 $$\lambda_1(x) \leq \dots \leq
\lambda_s(x) \leq 0 < \lambda_{s+1}(x)\leq \dots\leq \lambda_{d}(x)$$ are the Lyapunov exponents at $x$, then
\begin{equation}\label{e.entformula}
h_{\mu_f}(f)=\int_M\sum\limits_{i=s+1}^{d} \lambda_i(x)\,d\mu_f(x).
\end{equation}
We will refer to this last equality as the {\em entropy formula} for $\mu_f$. We will see in
Lemma~\ref{le.lyapunov} that in our situation $f$ has all its Lyapunov exponents positive with respect
to~$\mu_f$. Hence, by Oseledets  Theorem  and the integrability of the Jacobian of $f$ with respect to $\mu_f$,
we have that the integral in \eqref{e.entformula} is equal to the integral of the Jacobian of $f$ with respect
to the measure $\mu_f$; see Proposition \ref{pr.entropyf}.

One of the key results to prove our main result on continuity of the SRB entropy is the following theorem which
establishes the relation between the entropy of the original map and the entropy of the induced map with respect
to the appropriate measures.

\cmt\label{t.entropia}  
 If $F$ is an induced map for $f$ for which (i$_1$), (i$_2$) and (i$_3$) hold,
then
 $$h_{\mu_f}(f)=\frac{1}{\mu_f^*(M)}\,h_{\mu_F}(F).$$
 \fmt
The proof of this result will be given in Section \ref{Entropy
formulas}.

\subsection{Continuity of entropy}
Let $\cf$ be a family of $C^k$ maps, for some fixed $k\ge 2$, from a manifold $M$ into itself. Assume that we
may associate to each $f\in\cf$  an induced Markov map $ F_f\colon \Delta\rightarrow \Delta$ defined on a same
ball $\Delta\subset M$. Given $f\in\cf$, let $\cp_f$ denote the partition into domains of smoothness of $F_f$,
and $\tau_f\colon\cp_f\rightarrow \ZZ^+$ be its return time function. Let also $\mu_{F_f}$ be the absolutely
continuous $F_f$-invariant probability measure, $\mu_f^*$ the measure obtained from $\mu_{F_f}$ as in
\eqref{eq.relmu}, and $\mu_f$ its normalization.
 For notational simplicity we will denote   the Markov
induced map associated to $f$ by $F$ and its absolutely continuous invariant probability measure  by $\mu_F$.

One of the main goals of this work is to study the continuous
variation  of the metric entropy with respect to $\mu_f$ with the
map $f\in\cf$. In order to be able to implement our strategy we
assume that  the following uniformity conditions hold:
 \begin{enumerate}
 \item[(u$_1$)]\; \emph{$\tau_f$ varies continuously in the $L^1$ norm with $f\in\cf$.}
\item[(u$_2$)]\;  \emph{$\kappa$ and $K$ associated to $F_f$ as in (i$_2$) and (i$_3$) may be chosen uniformly
for $f\in\cf$\,.}
 \end{enumerate}
As we shall  see in Proposition~\ref{pr.continua}, these
uniformity conditions assure in particular that the (unique)
absolutely continuous probability measure $\mu_F$ invariant by the
map $F$ varies continuously (in the $L^1$ norm) with $f\in\cf$.

\cmt\label{t.cont} If $\cf$ is a  family of $C^k$ ($k\ge 2$)
 maps  from the manifold $M$ into itself for which (u$_1$) and (u$_2$) hold, then
the entropy
 $
  h_{\mu_f}(f)
 $
varies continuously with $f\in\cf$.
 \fmt

Next  we introduce a family of maps and present sufficient conditions for the validity of the assumptions of the
previous theorem. As we shall see these conditions are verified in the set of  maps introduced in \cite{V}.

\subsection{Non-uniformly expanding maps}
Let $f\colon M\to M$ be a $C^2$ local diffeomorphism in the whole manifold $M$ except possibly in a set of
critical points \( \cc \subset M\).
 We say that \( \cc \) is a {\em non-degenerate critical set} if the following conditions hold.
 The first one says that there are constants $B>0$ and $\be>0$ such that for every $x\in
 M\setminus\cc$ one has
\begin{itemize}
 \item[(c$_1$)]
\quad $\displaystyle{{\|Df(x)\|}\geq {B}\dist(x,\cc)^{\be}}$.
\end{itemize}
Moreover, we assume that the functions $  \log|\det Df| $ and $ \log \|Df^{-1}\| $ are \emph{locally Lipschitz}
at points $ x\in M \setminus \mathcal C $, with Lipschitz constant depending on $ \dist (x, \mathcal C)$: for
every $ x,y \in M\setminus \cc$ with $\dist(x,y)<\dist(x,\cc)/2$ we have
\begin{enumerate}
\item[(c$_2$)] \quad $\displaystyle{\left|\log\|Df(x)^{-1}\|- \log\|Df(y)^{-1}\|\:\right|\leq
\frac{B}{\dist(x,\cc)^{\be}}\dist(x,y)}$;
 \item[(c$_3$)]
\quad $\displaystyle{\left|\log|\det Df(x)^{-1}|- \log|\det Df(y)^{-1}|\:\right|\leq
\frac{B}{\dist(x,\cc)^{\be}}\dist(x,y)}$.
 \end{enumerate}

Note that the above conditions  give, for the particular case of
critical points of one-dimensional maps, the usual definition of a
non-degenerate critical point. From now on we assume that the
critical sets of the maps we will be considering are always
non-degenerate.

Given any $\delta>0$ and $x\in M\setminus\cc$, we define the {\em
$\delta$-truncated distance\/} from $x$ to $\cc$ as
\begin{equation*}
\label{e.truncate} \dist_\delta(x,\cc)=\left\{
\begin{array}{ll} 1, & \text{if } \dist(x,\cc)\ge\delta; \\
                  \dist(x,\cc), & \text{otherwise}.
\end{array}\right.
\end{equation*}
We say that $f$ is {\em non-uniformly expanding}  if the following two conditions hold:
\begin{itemize}
\item[(n$_1$)] {\em there is $\lambda>0$ such that for Lebesgue almost every $x\in M$
 \begin{equation*}\label{NUE} 
    \limsup_{n\to\infty}\frac{1}{n}\sum_{i=0}^{n-1}
    \log \|Df({f^{i}(x))}^{-1}\|<-\lambda;
\end{equation*}}
\item[(n$_2$)] {\em for every $\epsilon>0$ there exists $\delta>0$ such that for Lebesgue almost  every $x\in M$
\begin{equation*} \label{e.faraway1}
    \limsup_{n\to+\infty}
\frac{1}{n} \sum_{j=0}^{n-1}-\log \dist_\delta(f^j(x),\cc) \le\epsilon.
\end{equation*}}
\end{itemize}
We often refer to (n$_2$) by saying that orbits have {\em slow recurrence} to the critical set $\cc$. In the
case that $\cc$ is equal to the empty set we simply ignore the slow recurrence condition.

\cre It is worthy to be stressed that slow recurrence condition is not needed in all its strength for our
results. In fact, condition~(n$_2$) is needed just for distortion control reasons. As observed in \cite[Remark
1.3]{Al2},  it is enough to have it for some sufficiently small $\epsilon> 0$ and conveniently chosen $\delta
> 0$; see also \cite[Proposition 3.5]{Al2} and \cite[Remark 3.6]{Al2}. \fre

Condition  (n$_1$) implies that
 the
\emph{expansion time} function
\[
\mathcal E(x) = \min\left\{N\ge1\colon\frac{1}{n} \sum_{j=0}^{n-1} \log \|Df({f^{j}(x))^{-1}}\| \leq
-\frac\lambda2, \quad\text{for all $n\geq N$}\right\}
\]
is defined and finite Lebesgue almost everywhere in \( M \). The \emph{recurrence time}\index{time!recurrence}
function
\[
\mathcal R(x) = \min\left\{N\ge 1: \frac{1}{n} \sum_{j=0}^{n-1} -\log \dist_\delta(f^j(x),\cc) \leq
2\varepsilon, \quad\text{for all $n\geq N$}\right\},
\]
is also defined and finite Lebesgue almost everywhere in \( M \)  if the slow recurrence condition~(n$_2$)
holds.

 We think
of $\mathcal E(x) $ and $\mathcal R(x) $  as the time we need to wait before the exponential derivative growth
kicks~in. These depend on asymptotic statements and we have no a-priori knowledge about how fast these limits
are approached or with what degree of uniformity for different points \( x \). We define the {\em tail set (at
time $n$)}
\begin{equation}\label{tailset}
    \Gamma_n^f=\big\{x\in M \colon \mathcal E(x) > n \ \text{ or } \ \mathcal R(x) > n \big\}.
\end{equation}
This is the set of points which at time \( n \) have not yet achieved either the uniform exponential growth or
the slow recurrence given by conditions (n$_1$) and (n$_2$). If the critical set is empty, we simply ignore the
recurrence time function and consider only the expansion time function in the definition of $\Gamma_n^f$.

It is proved in \cite{ABV} that every $C^2$ non-uniformly expanding map $f$ admits some SRB measure. Moreover,
it follows from \cite[Lemma 5.6]{ABV} that if $f$ is transitive, then  it has a unique SRB measure $\mu_f$ which
is ergodic and absolutely continuous with respect to the Lebesgue measure, whose
basin covers a full Lebesgue measure subset of points in $M$. 

The results in \cite{Al2} show that  if  the decay of the Lebesgue measure of $\Gamma_n^f$ holds with some
uniformity in $f\in\cn$, then the SRB measure $\mu_f$ varies continuously in the $L^1$ norm with $f\in\cn$. Here
we deduce  the continuity of the SRB entropy in the same context.

\cmc \label{t.NUE2} Let $\cn$ be a set of $C^k$ ($k\ge2$) transitive non-uniformly expanding maps (with same
constants $\epsilon$, $\delta$ and $\lambda$). If there are $C>0$ and $\gamma>1$ such that $\leb(\Gamma_n^f)\le
Cn^{-\gamma}$, for all $f\in\cn$ and $n\ge1$, then the entropy $h_{\mu_f}(f)$ varies continuously with
$f\in\cn$.
 \fmc

This is a direct consequence of Theorem~\ref{t.cont} and the fact that,  by \cite[Proposition~5.3]{Al2}, maps in
a family $\cn$ as in the hypotheses of Corollary~\ref{t.NUE2} necessarily admit induced maps for which
uniformity conditions (u$_1$) and~(u$_2$) hold. The proof of \cite[Proposition~5.3]{Al2} uses ideas from
\cite{ALP}, where  piecewise expanding induced maps for non-uniformly expanding maps are constructed.
Transitivity is a useful ingredient for that construction.

\subsubsection{Viana maps}
Here we present an open class $\cv$  of transformations where the assumptions of Corollary~\ref{t.NUE2} hold.
This is an open set of maps from the cylinder into itself  constructed in \cite{V}. As pointed out in that
paper, the choice of the cylinder $S^1\times \RR$ as ambient space is rather arbitrary, and  the construction
extends easily to more general manifolds. In what follows we briefly describe the maps in the set $\cv$, and
refer the reader to \cite{V, Al, AV, AA,Bu}  for more details.

Let $a_0\in(1,2)$ be such that the critical point $x=0$ is pre-periodic under iteration by the quadratic map
$p(x)=a_0-x^2$, and let $b:S^1\rightarrow \RR$ be a Morse function, for instance, $b(t)=\sin(2\pi t)$. We take
$S^1=\RR/\ZZ$. For each $\alpha>0$, consider the map
 $$
f_{\al}:S^1\times \RR\rightarrow S^1\times \RR,\quad
f_\al(\theta,x)=(\hat{g}(\theta), \hat q(\theta,x)),
 $$
where $\hat{g}$ is the uniformly expanding map of the circle defined by $\hat{g}(\theta)=d\theta$ (mod $1$), for
some integer $d\ge 2$, and $\hat q(\theta,x)=a(\theta)-x^2$ with $a(\theta)=a_0+\al b(\theta)$. We  take $\cv$
as a small $C^3$ neighborhood of $f_\alpha$, for some (fixed) sufficiently small $\alpha>0$.
 Observe
that each $f\in\cv$ has a whole curve of critical points near $\{x=0\}$ for small enough $\alpha>0$. The
 $C^3$ topology is used in \cite{V} in order to simplify some technical points. In particular,
it is possible to prove $C^2$ proximity of the critical sets for $C^3$ nearby maps. We do believe that the
results in \cite{V} and the subsequent works for Viana maps still hold in the $C^2$ topology.

One can easily check that for $\al>0$ small enough there is an interval $I\subset (-2,2)$ such that $f_\al$
sends $S^1\times I$ into the interior of $S^1\times I$. Thus, any map $f$ close to $f_\al$ still has $S^1\times
I$ as a forward invariant region, and so it has an attractor inside this invariant region. The attractor is
precisely the set $ \Lambda=\cap_{n\geq 0}f^n(S^1 \times I)$.

It is proved in \cite{Al} that any $f \in \cv$ admits some absolutely continuous ergodic invariant probability
measure. Moreover, the results  in \cite{AV} show that these systems have a unique SRB measure whose basin
covers a full Lebesgue measure set of points in $S^1\times I$, and the densities of these SRB measures vary
continuously in the $L^1$ norm with the map.  To obtain the uniqueness of the SRB measure, they prove that $f$
is {\it topologically mixing}, in a strong sense: for every open set $A \subset S^1 \times I$ there is some $n =
n(A) \in \mathbb{Z}^+$ such that $f^n(A) = \Lambda.$ In particular, maps belonging to $\cv$ are transitive.

The non-uniform expansivity of Viana maps is proved in \cite{V}. Specific rates for the decay of the tail set
are known in this case: there exist constants $C, \gamma> 0$ (uniformly in the whole set $\cv$) such that
 $$m (\Gamma^f_n) \leq C \exp (-\gamma \sqrt{n}),\quad \text{for all $f\in\cv$ and $n\ge 1$};$$
 see \cite[Section 2.4]{V} and
\cite[Section 6.2]{AA} for details. Thus we may apply Corollary~\ref{t.NUE2} to the set of Viana maps and derive
the following consequence.

\cmc The SRB entropy of Viana maps varies continuously with $f\in \cv$. \fmc

%

Let us remark that $\cv$ is an open set in the space of $C^3$  transformations from the cylinder $S^1\times I$
into itself, where each $f\in\cv$ has a curve of critical points. The conclusion on the continuity of the SRB
entropy in this higher dimensional case is completely different from the above mentioned case of one-dimensional
quadratic maps.

\section{Statistical stability}

Let $\cf$ be a family os maps as in Theorem~\ref{t.cont}.  The main goal of this section is to prove that
$\mu_F$ varies continuously with $f\in\cn$.
 In the next lemma we give in particular
a proof that an absolutely continuous  invariant measure for a piecewise expanding Markov map exists.
 For the sake of notational
simplicity we shall write
 $$J_f(x)=|\det Df(x)|\qand
 J_F(x)=|\det DF(x)|.$$

\cpr\label{l.density}
 There is $C_0>0$ such that for each $f\in\cf $
there exists an $F$-invariant absolutely continuous probability measure $\mu_F= \rho_F m$ with $C_0^{-1}\le
\rho_F\le C_0$.

\fpr

\dem We start the proof of the result with the following claim: there exists $K_0>0$ such  that given  $f \in
\cf$ , $k\geq 1$, an inverse branch $G: \Delta \rightarrow G(\Delta)$ of $F^{-k}$, and measurable sets $A,
B\subset\Delta$, then
\begin{equation}\label{eqdistorcao}
K_0^{-1} \frac{m(A)}{m(B)}
 \leq \frac{m(G(A))}{m(G(B))}
 \leq K_0 \frac{m(A)}{m(B)}.
\end{equation}
Indeed, observe that
$$\frac{m(A)}{m(B)} = \frac{\int_{G(A)} J_{F^k} dm}{\int_{G(B)} J_{F^k} dm} $$
We use (i$_3$) and show that there is $K_1>0$ (uniformly choosen in $\cf$) such that
\begin{equation} \label{dige}
K_1^{-1} \le \frac{J_{F^k}(y)}{J_{F^k}(z)} \le K_1
 \end{equation}
for every $y$, $z$ on the image of $G$. For this purpose observe that
\begin{eqnarray*}
\log \frac{J_{F^k}(y)}{J_{F^k}(z)} & = &\sum_{i=0}^{k-1} \log \frac{J_F(F^i(y))}{J_F(F^i(z))}\\ & \leq &
\sum_{i=0}^{k-1} \left|\frac{J_F(F^i(y))}{J_F(F^i(z))} - 1 \right|\\ & \leq & K \sum_{i=1}^{k} \dist (F^i(y),
F^i(z)) \\
 &\leq & K \sum_{i=1}^{\infty} \kappa^{i} L,
\end{eqnarray*}
where $L$ is the diameter of $M.$ Observe that the last upper bound is uniform in $\cf$.

 Now we use (\ref{dige})
to prove (\ref{eqdistorcao}). Fixing $z \in G(\Delta)$, it comes out that
$$
\frac{\int_{G(A)} J_{F^k} dm}{\int_{G(B)} J_{F^k} dm} \leq K_1^2\frac{ J_{F^k}(z) m(G(A))}{ J_{F^k}(z) m(G(B))},
$$
and with the same argument we prove  the other inequality of (\ref{eqdistorcao}) with $K_1^2=K_0$.

Using the claim we will see that every accumulation point $\mu_F$ of the sequence
$$
\mu_n = \frac{1}{n}\sum_{i=0}^{n-1} F_{*}^i m
$$
is an $F$-invariant probability absolutely continuous with respect to $m$, with density $\rho_F$ bounded from
zero and from infinity. In order to prove it, take $B=\Delta$ and fix $C_0 = K_0m(\Delta)^{-1}$. Since $m
(F^{-k}(A))$ is the sum of the terms $m(G(A))$ over all inverse branches $G:\Delta \rightarrow G(\Delta)$ of
$F^k$, it follows from \eqref{eqdistorcao} that
 $$
 C_0^{-1} m(A) \leq m(F^{-k}(A)) \leq C_0 m(A).
 $$
 This implies that, for every $n$, the density $\rho_n=d\mu_n /d  m$  satisfies $C_0^{-1} \leq \rho_n \leq
C_0$, and the same holds for the density of the accumulation point  $\mu_F$. \cqd

 \cle \label{l.cauda}
 Given $\epsilon>0$, there are $N\geq 1$ and
$\de=\de(\epsilon,N)>0$
 such that for $f\in\cf$
 $$
 \|f-f_0\|_{C^k}<\de\quad\Rightarrow\quad
 m\{\tau_f>N\}<\varepsilon.
 $$
 \fle

 \dem For the sake of notational simplicity we denote $\tau_f$ by $\tau$
 and $\tau_{f_0}$ by $\tau_0$.
Take any $\e>0$ and take $N\geq 1$ in such a way that
 $
 \|\Chi_{\{\tau_{0}>N\}}\|_1<\vare/2,
 $
 where $\Chi_A$ denotes the indicator of a set $A$.
 We have
 \begin{align*}
 m\{\tau>N\}&=
 \big\|\Chi_{\{\tau>N\}}\big\|_1\\
 &=
 \big\|\Chi_{\{\tau>N\}} - \Chi_{\{\tau_{0}>N\}}+
 \Chi_{\{\tau_0>N\}}\big\|_1\\
  &\leq
\big\|\Chi_{\{\tau>N\}} - \Chi_{\{\tau_{0}>N\}}\|_1+
 \|\Chi_{\{\tau_0>N\}}\big\|_1
 \end{align*}
 and so, if we take $\de>0$ sufficiently small
 then, by (u$_1$), taking $\|f-f_0\|_{C^k}<\de$,
 the first term in the sum above can also be made smaller
 than $\vare/2$.
 \cqd

\cpr\label{pr.continua}  The measure $\mu_F$ varies continuously (in the $L^1$-norm) with $f\in \cf$. \fpr

\dem Let $f_n$ be any sequence in $\cf$ converging to $f_0$ in the $C^k$ topology. For each $n\ge 0$, consider
$F_n\colon\Delta\to\Delta$ the induced Markov map associated to $f_n$. Denote by $\rho_n$ the density of the
$F_n$-invariant absolutely continuous probability measure. Proposition~\ref{l.density} gives that the sequence
of densities $\rho_n$ is relatively compact in $L^\infty(\Delta, m)$, and so it has some accumulation point
$\rho_\infty$ with $\|\rho_\infty\|_\infty\le C_0$. With no loss of generality we assume that the full sequence
$\rho_n$ converges to $\rho_\infty$ in the $L^1$-norm. We need to prove that $\rho_\infty=\rho_0$. We will do
this by showing that
  $$
\int (\varphi\circ F_0)\rho_\infty dm= \int\varphi\rho_\infty dm
  $$
for every continuous $\varphi\colon \Delta\to\RR$, and use the fact that $F_0$ has a unique absolutely
continuous invariant probability measure.  Given any $\varphi\colon M\rightarrow\RR$ continuous we have
 $$
\int \f \rho_{n} dm\rightarrow \int \f  \rho_{\infty}dm\quad \mbox{when}\quad n\rightarrow\infty.
 $$ On
the other hand, since $\rho_{n}$ is the density of an $F_{n}$-invariant probability measure we have
 $$
 \int \f \rho_n dm=\int
(\f\circ F_{n})\rho_n dm\quad\mbox{for every }n\ge0.
 $$
 So, it suffices to prove that
 \begin{equation}
 \int (\f\circ F_{n})\rho_n dm
 \rightarrow \int (\f\circ F_0) \rho_\infty dm
 \quad\mbox{when}\quad n\rightarrow\infty.
 \end{equation}
 We have
 \begin{eqnarray*}
 \lefteqn{
 \big|\int (\f\circ F_{n})\rho_ndm -\int (\f\circ F_0) \rho_{\infty}dm| \leq }\\
& &
 \big|\int (\f\circ F_{n})\rho_{n}dm
 - \int (\f\circ F_0) \rho_ndm\big|
+\big|\int (\f\circ F_0) \rho_ndm
 - \int (\f\circ F_0) \rho_{\infty}dm\big|.
 \end{eqnarray*}
Since $\rho_n$ converges to $\rho_\infty$ in the $L^1$-norm and $\varphi\circ F_0$ is bounded in $\Delta$, we
easily deduce that the second term in the sum above is close to zero for large $n$.

The only thing we are left to prove is that the  first term in the sum above converges to 0 when $n$ tends to
$\infty$.  That term is equal to
\begin{equation*}
 \big|\int (\f\circ F_{n}-\f\circ F_0) \rho_ndm\big|.
 \end{equation*}
 Since $(\rho_n)_n$ is bounded in the $L^\infty$ norm by Proposition~\ref{l.density}, all we are left to show is
 \begin{equation}\label{eq.last}
 \int \big|\f\circ F_{n}-\f\circ F_0\big| dm\to
 0,\quad\text{when}\quad n\to \infty.
 \end{equation}
Take any $\vare>0$. For each $n\ge0$ let $\tau_n$ denote the return time function of $f_n$.  By
Lemma~\ref{l.cauda} there are $N\geq 1$ and $n_1\in\NN$ such that
 $$
 n\ge n_1\quad\Rightarrow\quad
 m(\{\tau_n>N\})<\varepsilon.
 $$
We write the integral in \eqref{eq.last}  as
 \begin{equation}\label{eq.somatau}
 \int_{\{\tau_n>N\}} \big|\f\circ F_{n}-\f\circ F_0\big|
 dm+\int_{\{\tau_n\le N\}} \big|\f\circ F_{n}-\f\circ F_0\big| dm.
 \end{equation}
 The first integral in \eqref{eq.somatau} is bounded by
 $2\vare\|\varphi\|_\infty$ for $n\ge n_1$. Let us now estimate
 the second integral in \eqref{eq.somatau}.
Define
 $$
 A_n=\big\{x\in\Delta\colon \tau_n(x)= \tau_0(x)\big\}.
 $$
 Since $\tau_n$ takes only integer values, we have by
 (u$_1$) that there is some $n_2\in\NN$ such that
  $$m(\Delta\setminus A_n)\le \vare,\quad\text{for each $n\ge n_2$}.
  $$
  Observe that for each $x\in A_n $ we have $
  F_n(x)=f_n^{\tau_0(x)}(x).$ Thus we may write
  $$
  \int_{\{\tau_n\le N\}} \big|\f\circ F_{n}-\f\circ F_0\big| dm
  \le \int_{\{\tau_0\le N\}} \big|\f\circ f_{n}^{\tau_0}-\f\circ f_0^{\tau_0}\big| dm
  +\int_{\Delta\setminus A_n} \big|\f\circ F_{n}-\f\circ F_0\big| dm.
  $$
Since $f_n\to f_0$ in the $C^k$ topology, there is $n_3\in \NN$ such that for $n\ge n_3$
 $$
 \int_{\{\tau_0\le N\}}\big|\f\circ f_{n}^{\tau_0}-\f\circ f_0^{\tau_0}\big| dm
 \le \vare m(\{\tau_0\le N\}).
 $$
 On the other hand, for $n\ge n_2$
 $$
 \int_{\Delta\setminus A_n} \big|\f\circ F_{n}-\f\circ F_0\big| dm
 \le 2\vare\|\varphi\|_\infty.
 $$
Thus we have for $n\ge \max\{n_1,n_2,n_3\}$
 $$
 \int \big|\f\circ F_{n}-\f\circ F_0\big| dm\le
 \vare\big(4\|\varphi\|_\infty+m(\{\tau_0\le N\})\big).
 $$
This proves \eqref{eq.last} since $\vare>0$ has been taken arbitrarily.
 \cqd

 \section{Entropy formulas}\label{Entropy formulas}

In this section we prove  Theorem~\ref{t.entropia}.  Let $F\colon\Delta\to\Delta$ be a piecewise expanding
Markov map and $\mu_F$ its absolutely continuous invariant probability measure. Since the Lypunov exponents of
the induced map $F$ (with respect to the measure $\mu_F$) are all positive, then  the next lemma shows, in
particular, that the Lyapunov exponents of $f$ (with respect to the measure $\mu_f$) are all positive.

%
%

\cle \label{le.lyapunov} If $\lambda$ is a Lyapunov exponent of $F$, then $\lambda/\bar\tau$ is a Lyapunov
exponent of $f$, where $\bar\tau=\int_\Delta\tau_f d\mu_F$.
 \fle

 \dem Let $n$ be a positive integer.  We have for each $x\in\Delta$
  $$F^n(x)=f^{S_n(x)}(x),\quad\text{where $S_n(x)=\sum_{i=0}^{n-1}\tau_f(F^i(x))$}.$$
 As $S_n(x) = S_n(y)$ for Lebesgue almost every $x\in \Delta$ and $y$ near enough $x$, we can
 take derivatives in the above equation and conclude that if $v\in T_xM$ then
  \begin{equation}\label{eq.ly}
  \frac{1}{S_n(x)}\log\|Df^{S_n(x)}(x)v\|=
  \frac{n}{nS_n(x)}\log\|DF^{n}(x)v\|.
  \end{equation}
Since $\mu_F$ is an ergodic measure, we have by Birkhoff's ergodic theorem
 \begin{equation}\label{eq.bi}
 \lim_{n\to\infty}\frac{S_n(x)}{n}=\int_\Delta\tau_fd\mu_F=\bar\tau
 \end{equation}
for Lebesgue almost every $x\in\Delta$ (recall that $\mu_F$ is equivalent to Lebesgue measure). Attending to
\eqref{eq.ly} and \eqref{eq.bi} the proof follows.
 \cqd

 \cpr\label{pr.entropyf}
The entropy formula holds for $\mu_f$, i.e.
 $\displaystyle h_{\mu_f}(f)=\int_M \log J_f\,d\mu_f. $
\fpr

\dem  As a consequence of Lemma \ref{le.lyapunov}, the fact that the Lyapunov exponents of $F$ with respect to
$\mu_F$ are all positive implies that all Lyapunov exponents of $f$ with respect to $\mu_f$ are also positive.
By the entropy formula
$$h_{\mu_f} = \int_M \sum\limits_{i=1}^d \lambda_i d \mu_f.$$  Now the integrability of $\log Jf$ with respect
to $\mu_f$ allows us to use Oseledets Theorem and rewrite the above equality as required in the above
proposition.
  \cqd

The proof of the next proposition uses fairly standard methods in ergodic theory.

 \cpr\label{pr.entropyF} If $F\colon \Delta\to\Delta$ is a piecewise expanding map for which (i$_1$), (i$_2$) and (i$_3$) hold, then
 $$h_{\mu_F}(F)=\int_\Delta \log J_F\,d\mu_F. $$
\fpr

\dem
First we observe that the measure $\mu_F$ is ergodic. We shall apply Shannon-McMillan-Breiman theorem for the
generating partition $\cp$ consisting of the smoothness domains of $F$. Take a generic point $x\in \Delta$. We have
\begin{equation}
  h_{\mu_F}(F)= h_{\mu_F}(F,\cp) = \lim_{n \rightarrow \infty} \frac{-1}{n} \log \mu_F(\cp_n(x)) = \lim_{n \rightarrow \infty} \frac{-1}{n} \log m(\cp_n(x)).
\end{equation}
The last equality comes from the fact that $m$ and $\mu_F$ are equivalent measures with uniformly bounded
densities. Now observe that each $\cp_n(x) $ is equal to some $ G(\Delta)$, where $G$ is an inverse branch of
$F^n$. Hence we have

\begin{equation}
 m(\Delta) = \int_{G(\Delta)} J_{F^n} dm.
\end{equation}
By the distortion estimate obtained in the proof of the Proposition \ref{l.density} we conclude that
$$
 K_1^{-1} \leq  m(G(\Delta)) J_{F^n} (x) \leq K_1.
$$
 By the above inequality we deduce that
 $$
  \lim_{n \rightarrow \infty} \frac{-1}{n} \log m(\cp_n(x))= \lim_{n \rightarrow \infty} \frac{1}{n} \log J_{F^n}(x) =
  \lim_{n \rightarrow \infty} \frac{1}{n} \sum_{i=0}^{n-1}\log J_F(F^i(x)) = \int \log J_F d\mu_F,
 $$
 where the last equality holds by Birkhoff ergodic theorem.
\cqd


Now we give a lemma with the aid of which we shall prove Theorem~\ref{t.entropia}. \cle\label{le.igual} If $F$
is an induced piecewise expanding Markovian map for $f$, then
$$\int_\Delta\log J_F\,d\mu_F=\int_M\log J_f\,d\mu_f^*.$$

\fle

\dem We define for each $n\ge 1$
 $$P_n=\{\omega\in\cp\colon \tau(\omega)=n\}.$$
Observe that for each $x\in P_n$ we have $F=f^n$. So, by the chain rule,
 $$J_F(x)=J_f(f^{n-1}(x))\cdots J_f(f(x))\cdot J_f(x).$$
Thus we have for each $n\ge 1$
 \begin{eqnarray*}
 \lefteqn{\int_{P_n}\log J_F d\mu_F= }\hspace{0cm}\\
 &=&\int_{P_n}\log  J_f\circ f^{n-1}d\mu_F+\cdots+\int_{P_n}\log
 J_f\circ f d\mu_F + \int_{P_n}\log
 J_f d\mu_F\\
 &=& \int_M\log  J_f\,d\left(f^{n-1}_*(\mu_F\vert P_n)\right)+\cdots + \int_M\log
 J_f \,d\left(f_*(\mu_F\vert P_n)\right)+\int_M\log
 J_f\,d(\mu_F\vert P_n).\\
 \end{eqnarray*}
 Using this we deduce
 \begin{eqnarray*}
 \int_\Delta\log J_F d\mu_F &= &\sum_{n=1}^\infty\int_{P_n}\log J_F\,d\mu_F\\
  &=&
 \sum_{n=1}^\infty\sum_{j=0}^{n-1}\int_M\log J_f\,d\left(f^{j}_*(\mu_F\vert P_{n})\right)\\
  &=&
 \sum_{n=0}^\infty\int_M\log  J_f\,d\left(f^{n}_*(\mu_F\vert \{\tau>n\})\right)\\
  &=&
  \int_M\log  J_f\,d\left(\sum_{n=0}^\infty f^{n}_*(\mu_F\vert
\{\tau>n\})\right).
 \end{eqnarray*}
By \eqref{eq.relmu} we have
$$
\int_M\log  J_f\,d\left(\sum_{n=0}^\infty f^{n}_*(\mu_F\vert \{\tau>n\})\right)= \int_M\log J_f d\mu_f^*,
$$
and so we have proved the result.
 \cqd

 Since the
entropy formula holds for $\mu_f$ by Proposition~\ref{pr.entropyf}, then using Proposition~\ref{pr.entropyF} and
Lemma~\ref{le.igual} we obtain
 \begin{eqnarray}
 h_{\mu_f}(f) &=& \int_M\log J_f\,d\mu_f\nonumber\\
 &=&
 \frac{1}{\mu_f^*(M)}\,\int_M\log J_f\,d\mu_f^*\nonumber\\
 &=&
 \frac{1}{\mu_f^*(M)}\,\int_\Delta\log J_F\,d\mu_F\label{igual}\\
 &=&
 \frac{1}{\mu_f^*(M)} h_{\mu_F}(F).\nonumber
 \end{eqnarray}
 This proves  Theorem~\ref{t.entropia}.

\section{Continuity of entropy}

In this section we prove Theorem ~\ref{t.cont}. Let $\cf$ be a
family of $C^k$ maps, $k\ge 2$, from the manifold $M$ into itself
for which (u$_1$) and (u$_2$) hold. We are implicitly assuming
that we have some $\Delta\subset M$ and, associated to each
$f\in\cf$, a piecewise expanding Markov induced   map
$F\colon\Delta\to\Delta$. By \eqref{igual}, in order to prove
Theorem~\ref{t.cont},  we just have to show that both $\mu_f^*(M)$
and $\int_\Delta\log J_F\,d\mu_F$ vary continuously with
$f\in\cf$.

Take an arbitrary $f_0\in\cf$ and let $f_n$ be any sequence in
$\cf$ converging to $f_0$ in the $C^k$ topology.  For each $n\ge
0$, let $F_n\colon\Delta\to\Delta$ be the induced map associated
to $f_n$, and let $\tau_n\colon \Delta\to\NN$ be the respective
return time function. Denote by $\rho_n$ the density of the
absolutely continuous $F_n$-invariant probability measure
$\mu_{F_n}$. Consider also for $n\ge 0$ the absolutely continuous
$f_n$-invariant measure $\mu_n^*$ obtained as in \eqref{eq.relmu}
from $\mu_{F_n}$:
\begin{equation*}\label{eq.relmu2}
 \mu_n^*
=\sum_{j=0}^{\infty}f_{\ast}^j\left(\mu_{F_n}\mid \{\tau_f>j\}\right).
 \end{equation*}
The continuous variation of $\mu_f^*(M)$ and $\int_\Delta\log
J_F\,d\mu_F$ with $f\in\cf$ will follow from
Proposition~\ref{jadid} and Proposition~\ref{contlog} below. We
start with an abstract lemma.

\cle\label{le.preliminar} Let $(\varphi_n)_n$ be a  bounded sequence in $L^\infty(m)$. If $\varphi_n\to\varphi$
in the $L^1(m)$ norm and $\psi\in L^1(m)$, then
 $$\int\psi(\varphi_n-\varphi)dm\to0,\quad\text{when $n\to\infty$}.$$ \fle

\dem Take any $\epsilon > 0$. Let $C>0$ be an upper bound for $\|\varphi_n\|_\infty$.
 Since $\psi\in L^1(m)$, there is $\delta>0$ such that for
any Borel set $B\subset M$
 \begin{equation}\label{eqz}
 m(B)<\delta\quad\Rightarrow\quad\int_B|\psi|dm<\frac\epsilon{4C}.
  \end{equation}
 Define for each $n\ge 1$
 $$
 B_n= \left\{x\in \Delta\colon |\varphi_n(x)-\varphi_{0}(x)|> \frac\epsilon{2\|\psi\|_1} \right\}.
 $$
Since $\|\varphi_n-\varphi_{0}\|_1 \rightarrow 0$ when
$n\to\infty$, then there is $n_0\in\NN$ such that $m ( B_n)
<\delta$ for every $n\ge n_0$. Taking into account  the definition
of $B_n$, we may write
\begin{align*}
\int |\psi||\varphi_n-\varphi_{0}|dm &= \int_{B_n}
|\psi||\varphi_n-\varphi_{0}|dm + \int_{\Delta\setminus B_n}
|\psi||\varphi_n-\varphi_{0}|dm\\
 &\leq 2C \int_{B_n}
|\psi| dm + \frac\epsilon{2\|\psi\|_1} \int_{\Delta\setminus B_n}
|\psi|dm.
\end{align*}
Then, using \eqref{eqz}, this last sum is upper bounded by $\epsilon$, as long as $n\ge n_0$. 
\cqd

\cpr  \label{jadid} $\mu_n^*(M)$ converges to $\mu_0^*(M)$ when
$n\to\infty$. \fpr

\dem  Recall that we have for every $n\ge0$
 \begin{equation*}
 \mu_n^*(M)
=\sum_{j=0}^{\infty}\mu_{F_n}\left( \{\tau_n>j\}\right)= \int
\tau_n d\mu_{F_n}.
 \end{equation*}
Hence
 $$
 |\mu_{n}^*(M)-\mu_{0}^*(M)| \leq \int
 |\tau_n\rho_n-\tau_{0}\rho_{0}|dm.
 $$
Now we write
 \begin{equation}\label{eq.z}
    \int |\tau_n\rho_n-\tau_{0}\rho_{0}|dm \leq \int|\tau_{0}||\rho_n-\rho_{0}|dm
    +\int |\tau_n-\tau_{0}||\rho_{n}| dm.
  \end{equation}
Let us first control the first term on the right hand side of
~\eqref{eq.z}. If we take $\psi=\tau_0$ and $\varphi_n=\rho_n$ for
each $n\ge 0,$  then, by Proposition~\ref{l.density} and
Proposition \ref{pr.continua}, these functions are in the
conditions of Lemma~\ref{le.preliminar}. Hence, the first term on
the right hand side of \eqref{eq.z} converges to 0 when
$n\to\infty$. We just have to notice that
 \begin{equation}\label{left}
 \int |\tau_n-\tau_{0}||\rho_{n}| dm\to0,\quad\text{when}\quad n\to\infty.
 \end{equation}
In fact, since $(\rho_n)_n$ is uniformly bounded by Proposition~\ref{l.density}, then  hypothesis (u$_1$)
assures that \eqref{left} holds.
  \cqd

At this point we have proved the continuous variation of $\mu_f(M)$ with $f\in\cf$, thus attaining the first
step in the proof of Theorem~\ref{t.cont}. The next step is to prove the continuous variation of
$\int_\Delta\log J_F\,d\mu_F$ with $f\in\cf$. We start with an auxiliary lemma.

\cle\label{le.majora} There is $C>0$ such that $\log J_{F_n}\le C\tau_n$ for every $n\ge 0$. \fle

\dem Define   $K_n=\max_{x\in M}\{J_{f_n}(x)\}$, for each $n\ge 0$.
  By the compactness  of $M$ and the
continuity on the first order derivative,  there is $K > 1$ such that $K_n\le K$ for all $n\ge 0$. We have
 $$
 J_{F_n}(x)=\prod_{j=0}^{\tau_n(x)-1}J_{f_n}(f_n^j(x))\le K^{\tau_n(x)}.
 $$
Hence
 $$0 < \log J_{F_n}(x) \le \tau_n(x)\log K.$$
We just have to take $C=e^K$. \cqd

The previous result gives in particular the integrability of $\log J_f $ with respect to Lebesgue measure, under
the assumption of the integrability of  $\tau_f$.  In the proof of the next proposition we also obtain the
continuous variation  of $\log J_f $ in the $L^1(m)$ norm with $f\in\cf$, as explicitly stated  in
\eqref{contlf} below.

 \cpr\label{contlog}
$\displaystyle{\int\log J_{F_n}\,d\mu_{F_n}}$ converges to $\displaystyle{\int\log J_{F_0}\,d\mu_{F_0}}$, when
$n\to\infty$.
 \fpr

 \dem
First we write
 \begin{eqnarray*}
 \lefteqn{
\big| \int \log J_{F_0} d \mu_{F_0} -  \int \log
J_{F_n} d \mu_{F_n} \big| \leq }\\
& &\hspace{4cm}
 \big| \int{(\log J_{F_n} - \log J_{F_0} ) \rho_n} dm \big| +
\big| \int{ (\rho_n -\rho_0)\log J_{F_0}} dm \big|.
 \end{eqnarray*}
It follows from Proposition~\ref{l.density}, Proposition~\ref{pr.continua} and Lemma~\ref{le.majora}  that if we
take $\varphi_n=\rho_n$ and $\psi=\log J_{F_0}$ then these functions are in the conditions of
Lemma~\ref{le.preliminar}. Thus, it is enough to show that
 \begin{equation}\label{contlf}
  \int{\big|\log J_{F_n} - \log J_{F_0} \big|} dm
 \to0,\quad\text{when $n\to\infty$.}
 \end{equation}
Take any $\vare>0$.  Since $\tau_0\in L^1(m)$, there is $N\geq 1$  such that
 \begin{equation}\label{eq.tau1}
 \int_{\{\tau_0>N\}}\tau_0 dm<\varepsilon.
  \end{equation}
 We then write
  \begin{eqnarray}\label{eq.duas}
 \lefteqn{\int|\log J_{F_n} - \log J_{F_0} | dm \,= }\\
& &\hspace{1cm}  \int_{
 \{\tau_n>N\}}|\log J_{F_n} - \log J_{F_0} | dm +
  \int_{\{\tau_n\le N\}}|\log J_{F_n} - \log J_{F_0} | dm .
  \end{eqnarray}
Let us start by controlling  the first in this last sum.  Using Lemma~\ref{le.majora} we obtain
 \begin{equation}\label{eq.tres} \int_{
 \{\tau_n>N\}}|\log J_{F_n} - \log J_{F_0} | dm
 \le
  \int_{
 \{\tau_n>N\}}\tau_n dm+\int_{
 \{\tau_n>N\}} \tau_0 dm.
 \end{equation}
One has
\begin{equation}\label{eq.tau2}
\Chi_{\{\tau_n>N\}}\tau_n\le \Chi_{\{\tau_0>N\}}\tau_0+
|\Chi_{\{\tau_n>N\}}-\Chi_{\{\tau_0>N\}}|\tau_0+\Chi_{\{\tau_n>N\}}
 |\tau_n-\tau_0|.
\end{equation}
 Choosing $n$ sufficiently large, we have
\begin{equation}\label{eq.tau3}
\int\Chi_{\{\tau_n>N\}}|\tau_n-\tau_0|dm\le \int|\tau_n-\tau_0|dm <\vare.
\end{equation}
On the other hand, applying Lemma~\ref{le.preliminar} to $\varphi_n=\Chi_{\{\tau_n>N\}}$, for $n\ge 0$, and
$\psi=\tau_0$ we also have for large $n$
\begin{equation}\label{eq.tau4}
\int|\Chi_{\{\tau_n>N\}}-\Chi_{\{\tau_0>N\}}|\,\tau_0 dm <\vare.
\end{equation}
It follows from \eqref{eq.tau1}, \eqref{eq.tau2}, \eqref{eq.tau3} and \eqref{eq.tau4} that for large $n$
 \begin{equation}\label{eq.con1}\int_{
 \{\tau_n>N\}}\tau_n dm<3\vare.\end{equation}
 Also from \eqref{eq.tau1} and \eqref{eq.tau4}
 \begin{equation}\label{eq.con2}
 \int_{
 \{\tau_n>N\}} \tau_0 dm \le \int|\Chi_{\{\tau_n>N\}}-\Chi_{\{\tau_0>N\}}|\,\tau_0 dm+\int_{\{\tau_0>N\}}\tau_0 dm
<2\vare.\end{equation} Hence, from \eqref{eq.tres}, \eqref{eq.con1} and \eqref{eq.con2} we deduce that for large
$n$
\begin{equation}\label{eq.tres2} \int_{
 \{\tau_n>N\}}|\log J_{F_n} - \log J_F | dm
 <5\vare.
 \end{equation}

Let us now estimate the second term in \eqref{eq.duas}. Letting $C>0$ be the constant given by
Lemma~\ref{le.majora}, take $\delta>0$ such that
 \begin{equation}\label{eq.rad1} \int_B C(N +\tau_0) dm<\vare,\quad\text{whenever $m(B)<\delta$}.
 \end{equation}
 For each $n\in \NN$ define
 $$
 A_n=\big\{x\in\Delta\colon \tau_n(x)= \tau_0(x)\big\}.
 $$
 Since $\tau_n$ takes only integer values, we have by
 (u$_1$)
  \begin{equation}\label{eq.rad2}
  m(\Delta\setminus A_n)\le \delta, \quad\text{for large $n$.}
 \end{equation}
  Observe that for each $x\in A_n $ we have $
  F_n(x)=f_n^{\tau_0(x)}(x).$ Thus we may write
  \begin{eqnarray*}
  \lefteqn{\int_{\{\tau_n\le N\}} |\log J_{F_n} - \log J_{F_0} | dm \le }\\
  & &\hspace{2cm}
  \int_{A_n\cap \{\tau_n \le N\}} |\log J_{ f_{n}^{\tau_0}}-\log J_{f_0^{\tau_0}}| dm
  +\int_{\{\tau_n\le N\}\setminus A_n} |\log J_{F_{n}}-\log J_{F_0}| dm.
  \end{eqnarray*}
Note that  by (i$_2$) we have  $J_{ f_{n}^{\tau_0}}\ge 1$ for every $n\ge 0$. Hence, the first integral in the
last sum can be made arbitrarily small if we take $n$ sufficiently large. On the other hand, we have by
Lemma~\ref{le.majora}
 $$
 \int_{\{\tau_n\le N\}\setminus A_n} |\log J_{F_{n}}-\log J_{F_0}| dm \le
 \int_{\Delta\setminus A_n} C(N +\tau_0) dm
 $$
It follows from \eqref{eq.rad1} and \eqref{eq.rad2} that this last quantity can be made smaller than $\vare>0$,
as long as $n$ is take sufficiently large.
 \cqd

%

\end{document}